\documentclass[11pt]{article}
\usepackage{amssymb, amsmath}
\title{ Remarks on the Extremal Functions for the
Moser-Trudinger Inequalities}
\author{Yuxiang Li
\\ {\it\small Department of Mathematical Science ,Tsinghua University,
Beijing, P.R.China, 100084}
\\ {\it\small E-mail address:
yxli@math.tsinghua.edu.cn}}
\date{}
\begin{document}
\maketitle
\begin{abstract}
We will show in this paper that if $\lambda$ is very close to 1,
then
$$I(M,\lambda,m)=
\sup_{u\in H^{1,n}_0(M) ,\int_M|\nabla u|^ndV=1}\int_\Omega
(e^{\alpha_n
|u|^\frac{n}{n-1}}-\lambda\sum\limits_{k=1}^m\frac{|\alpha_nu^\frac{n}{n-1}|^k}
{k!})dV,$$ can be attained, where $M$ is a compact manifold with
boundary. This result gives a counter example to the conjecture of
de Figueiredo, do \'o, and Ruf in their paper titled "On a
inequality by N.Trudinger and J.Moser and related elliptic
equations" (Comm. Pure. Appl. Math.,{\bf 55}:135-152, 2002).
\newline
\\{\bf Keywords:} Moser-Trudinger inequality, extremal function.
\newline
\\{\bf MR (2000) Subject Classification:} 58J05
\end{abstract}

\section{Introduction:}
Recall the Sobolev embedding: $H^{1,p}(\Omega)\hookrightarrow
L^\frac{np}{n-p}(\Omega)$ where $\Omega$ is a bounded domain in
$\mathbb{R}^n$ and $p\in (1,n)$.

When $p=n$, the critical growth is given by
the Moser-Trudinger inequality ([5]) which can be expressed as:
$$\sup_{u\in H^{1,n}_0(\Omega),||\nabla u||_{{L^n}(\Omega)=1}}\int_\Omega(e^{\alpha
|u|^\frac{n}{n-1}}-1)dx=c_n(\alpha)|\Omega|,\eqno (1.2)$$
then, denoting by $\omega_{n-1}$ the measure of unit
sphere in $\mathbb{R}^n$, one has
$$c_n(\alpha)<+\infty\;\; for\;\; 0<\alpha\leq\alpha_n=n\omega_{n-1}^\frac{1}{n-1},$$
$$c_n(\alpha)=+\infty\;\; for\;\; \alpha>\alpha_n.$$

Unlike Sobolev inequalities, the extremal functions for (1.2) with
$\alpha=\alpha_n$ exist generally. The first result for the
existence of extremal functions belongs to Carleson and Chang
([1]) who proved that (1.2) is attained when $\Omega=B_1(0)$ and
$\alpha=\alpha_n$. Their work has been extended to bounded domains
by Flucher([3]) and K-C, Lin([8]).

In [2] de Figueiredo et. al. gave a new proof for Carleson-Chang's
theorem. In fact, they got a generalized result which states that
$$\sup_{u\in H^{1,n}_0(B_1(0)),||\nabla u||_{L^n(B_1(0))}=1}
\int_{B_1(0)}(e^{\alpha_n
|u|^\frac{n}{n-1}}-\lambda|u|^\frac{n}{n-1})dx$$
is attained for any $\lambda<\alpha_n$.
Moreover, they proved the following:
\\{\bf Theorem A} {\it let $F(t)=F(|t|)$ be a function increasing on $\mathbb{R}^{+}$
which satisfies
$$1\leq F(t)\leq e^{\alpha_n|t|^\frac{n}{n-1}},\;and\;
\lim_{t\rightarrow+\infty}\frac{F(t)}{e^{\alpha_n|t|^\frac{n}{n-1}}}
=1,$$
then either
$$\sup_{u\in H^{1,n}_0({B_1(0)}),
\int_{B_1(0)}|\nabla u|^n=1}\int_{B_1(0)}
F(u)dx$$
can be attained, or
$$\sup_{u\in H^{1,n}_0({B_1(0)}),
\int_{B_1(0)}|\nabla u|^n=1}\int_{B_1(0)}
F(u)dx=e^{1+1/2+\cdots+1/(n-1)}+|B_1(0)|.\eqno (1.3)$$}

de Figueiredo et. al. then raised the following open problem in [1]:

{\bf Open problem: }{\it Show that
$$\sup_{\int_{B_1}|\nabla u|^n=1,u\in H^{1,2}_0(B_1)}\int_{B_1}F(u)dV$$
is not attained for $F(t)$ of the form
$$F(t)=e^{\alpha_n |t|^\frac{n}{n-1}}-g(t),$$
with
$$\lim_{t\rightarrow \infty}
\frac{g(t)}{e^{\alpha_n|t|^\frac{n}{n-1}}}=0\;\;and\;\;
g(t)\geq \alpha_n|t|^{\frac{n}{n-1}}.$$}

However, this problem is much more complicated than they expected,
we will show in this paper that the open problem is not true when
$g(t)$ is very close to
$$g_m(t)=\sum\limits_{k=1}^m\frac{|\alpha_nt^\frac{n}{n-1}|^k}
{k!}.$$

Let
$$F_{\lambda,m}(t)=e^{\alpha_n|t|^\frac{n}{n-1}}-
\lambda g_m(t),$$
and
$$I(\Omega,\lambda,m)=
\sup_{u\in H^{1,n}_0(\Omega) ,\int_\Omega|\nabla u|dV=1}\int_\Omega F_{\lambda,m}
(u)dV,$$
we shall prove that $I(\Omega,\lambda,m)$ is attained
when $\lambda$ is close to 1, in fact we shall prove it is attained
even when $\Omega$ is a manifold.

The author have adopted blow up analysis to study a similar problem. In
[7] the author
extended Carleson-Chang's result to
compact Riemannian manifolds. With the same
method in [7], we can also get a
generalization of Theorem A as follows:
\vspace{0.7ex}
\\{\bf Theorem 1.1:} {\it Let $(M,g)$ be a compact Riemannian manifold with
boundary. If
$I(M,\lambda,m)$ can not be attained, then we have
$$I(M,\lambda,m)=\mu(M)+\frac{\omega_{n-1}}{n}
e^{\alpha_nS_p+1+1/2+\cdots+1/(n-1)},$$
for some $p\overline{\in}\partial M$, where $\mu(M)$ is the measure of M.}

Here we should explain what $S_p$ is.
Let $G_p$ be the n-Lapalace Green function which is defined by
$$\left\{\begin{array}{l}
           -div(|\nabla G_p|^{n-2}\nabla G_p)=\delta_p\\
          {G_p}|_{\partial M}=0.
          \end{array}\right.\eqno (1.4)$$
By the asymptotic expansion theorem in [4] (the author has
shown in [7] that their proof still works on manifolds), the function
$G_p(x)+\frac{1}{\alpha_n}\log{dist^n(x,p)}$, is continuous at p.
$S_p$ is just its value at p, i.e.
$$S_p=\lim_{x\rightarrow p}(G_p(x)+\frac{1}{\alpha_n}\log{dist^n(x,p)}).\eqno (1.5)$$

Since the proof for Theorem 1.1 is the same as
the one done in [7] (also in [6]), we only give an outline of
the proof from which one can see the idea why we construct the
function sequence (3.5).

Then we will show that
\vspace{0.7ex}
\\{\bf Theorem 1.2:} {\it $I(M,1,m)>\mu(M)+\frac{\omega_{n-1}}{n}
e^{\alpha_nS_p+1+1/2+\cdots+1/(n-1)}$.}
\vspace{0.7ex}

Clearly, the function $f(\lambda)=I(M,\lambda,m)$
is continuous for fixed integer $m$, therefore we have
\vspace{0.7ex}
\\{\bf Main Theorem:} {\it There is a constant $\lambda_0>1$, s.t.
$I(M,\lambda,m)$ can be attained on $[0,\lambda_0)$.}

\section{The outline of the proof for Theorem 1.1}
Let $u_k\in H^{1,n}_0(M)$ which satisfies $u_k\geq 0$, $\int_M|\nabla u_k|^ndV_g=1$ and
$$\int_M(e^{\beta_k|u_k|^\frac{n}{n-1}}-\lambda g_m(u_k))
dV_g=\sup_{v\in H^{1,n}_0(M),\int_M|\nabla v|^ndV_g=1}
  \int_M(e^{\beta_k|v|^\frac{n}{n-1}}-\lambda g_m(v))dV_g,$$
where $\{\beta_k\}$ is an increasing sequence which converges to
$\alpha_n$. Then we have
$$-div|\nabla u_k|^{n-2}\nabla u_k=\frac{u_k^\frac{1}{n-1}
e^{\beta_ku_k^\frac{n}{n-1}}-\lambda g'_m(u_k)}{\lambda_k},$$
where
$$\lambda_k=\int_M(u_k^\frac{n}{n-1}e^{\beta_ku_k^\frac{n}{n-1}}-\lambda g_m'u_k)dV_g.$$
Let $c_k=\max\limits_{M}u_k(x)$, and $x_k\rightarrow p$. Then
using the same method in [7], we have $p\overline{\in}\partial M$.

Take a normal coordinate system $(U,x)$ around p.
Let $r_k^n=\frac{\lambda_k}{c_k^\frac{n}{n-1}e^{\beta_kc_k^\frac{n}
{n-1}}}$, then we have $r_k\rightarrow 0$ and
$$\frac{n}{n-1}\beta_kc_k^\frac{1}{n-1}(u_k(x_k+r_kx)-c_k)\rightarrow
-n\log(1+c_nr^\frac{n}{n-1})\eqno (2.1)$$
on any $B_L(0)$, where $c_n=(\frac{\omega_{n-1}}{n})^\frac{1}{n-1}$.

Moreover, we can get
$$c_k^\frac{1}{n-1}u_k\rightarrow G_p \eqno (2.2)$$
on any $M\setminus B_\delta(p)$ where $G_p$
is the Green function  defined by (1.4).

It follows from (2.1) that
$$\sup_{v\in H^{1,2}_0(M),\int_M|\nabla v|^ndV_g=1}F_{\lambda,m}(v)dV_g
=\lim_{k\rightarrow +\infty}
\int_M(e^{\beta_k|u_k|^\frac{n}{n-1}}-\lambda g_m(u_k))
dV_g=\mu(M)+
\lim_{k\rightarrow+\infty}\frac{\lambda_k}{c_k^\frac{n}{n-1}}.
$$
Then Theorem 1.1 can be obtained from an estimate for
$\frac{\lambda_k}{c_k^\frac{n}{n-1}}$.

\section{The proof of Theorem 1.2 :}
The following lemma will play an important role in
our computation.
\\{\bf Lemma 3.1} {\it Let $A_t=\{x\in M: {G_p}(x)>t\}$, then
as $t\rightarrow +\infty$, we have
$$\int_{\partial A_t}\frac{1}{|\nabla {G_p}|}\geq \omega_{n-1}^\frac{n}{n-1}
e^{-\alpha_nt+\alpha_nS_p}(1+O(e^{-\frac{2}{n}\alpha_nt})).$$}
\\{\bf Proof:} We have
$$-\frac{d|A_t|}{dt}=\int_{\partial A_t}\frac{dS_t}{|\nabla {G_p}|}.\eqno (3.1)$$
By H\"older inequality, we get
$$(\int_{\partial A_t}\frac{dS_t}{|\nabla {G_p}|})^{n-1}\int_{
\partial A_t}|\nabla {G_p}|^{n-1}dS_t\geq |\partial A_t|^n.\eqno (3.2)$$
Since $\nabla {G_p}=\frac{\partial {G_p}}{\partial n}$ on $\partial A_t$,
we have
$$\int_{\partial A_t}|\nabla {G_p}|^{n-1}dS_t=\int_{\partial A_t}
  |\nabla {G_p}|^{n-2}\frac{\partial {G_p}}{\partial n}\times 1dS_t=1.$$

Choose a normal coordinate system around p, then
we are able to assume that $g=\sum\limits_i(dx^i)^2+O(r^2)$, where
$r^2=dist^2(x,p)=\sum\limits_k(x^k)^2$.

Set
$$G_p(x)=-\frac{1}{\alpha_n}\log{r^n}+S_p+H(x),$$
where $H$ is a continuous function with $H(0)=0$.
Then for a point $x\in \partial A_t$, we have
$t=-\frac{1}{\alpha_n}\log{r^n(x)}+S_p+H(x)$, i.e.
$$r(x)=e^{-\frac{\alpha_n}{n}t+\frac{\alpha_nS_p+H(x)}{n}}.\eqno(3.3)$$
We can easily see that
$$|A_t|=\int_{A_t}\sqrt{|g|}dx=
\int_{A_t}(1+O(r^2))dx=(1+O(e^{-\frac{2\alpha_n}{n}t}))L(A_t),$$
where $L(\cdot)$ is the standard measure in $\mathbb{R}^n$. In the same way,
we have
$$|\partial A_t|=(1+O(e^{-\frac{2\alpha_n}{n}t}))L(\partial A_t).$$
Applying isopermetric inequality on $\mathbb{R}^n$, we get
$$|\partial A_t|^\frac{n}{n-1}\geq\alpha_n(1+O(e^{-\frac{2\alpha_n}{n}t}))|A_t|.$$
Thus, by (3.1), (3.2) we have
$$-\frac{d|A_t|}{dt}\geq\alpha_n(1+O(e^{-\frac{2}{n}\alpha_nt}))|A_t|.\eqno(3.4)$$
Hence, we have
$$\frac{d(e^{\alpha_nt}|A_t|)}{dt}\leq O(e^{-\frac{2}{n}\alpha_nt})(e^{
\alpha_nt}|A_t|)=O(e^{-\frac{2}{n}\alpha_nt}).$$
It is easy to see that
$\lim\limits_{t\rightarrow\infty}e^{\alpha_nt}|A_t|=
\frac{\omega_{n-1}}{n}e^{\alpha_nS_p}$, then we get
$$e^{\alpha_nt}|A_t|\geq\frac{\omega_{n-1}}{n}e^{\alpha_nS_p}+
O(e^{-\frac{2}{n}\alpha_nt}).$$
So, by (3.4) we have
$$\begin{array}{lll}
   -\frac{d|A_t|}{dt}&\geq&\alpha_n\frac{
    \omega_{n-1}}{n}e^{-\alpha_nt+\alpha_nS_p}(1+O(e^{-\frac{2}{n}\alpha_nt}))(1+
   O(e^{-\frac{2}{n}\alpha_nt}))\\
  &=&\omega_{n-1}^\frac{n}{n-1}e^{-\alpha_nt+\alpha_nS_p}(1+O(e^{-\frac{2}{n}\alpha_nt})).
  \end{array}$$
\begin{flushright}
{\small $\Box$}
\end{flushright}

Now, we are able to construct a function sequence $\{u_\epsilon\}\subset
H^{1,n}_0(M)$ which satisfies
$$\int_MF_{1,m}(u_\epsilon)dV_g>
\mu(M)+\frac{\omega_{n-1}}{n}e^{\alpha_nS_p+1+1/2+
 \cdots+/1(n-1)}.$$
Let
$$f_\epsilon(t)=\left\{\begin{array}{ll}
                 {C_\epsilon}-\frac{(n-1)\log(1+c_n\frac{e^{-\frac{\alpha_n}{n-1}t}
                    }{\epsilon^\frac{n}{n-1}})+{\Lambda_\epsilon}}{\alpha_n
                   {C_\epsilon}^\frac{1}{n-1}}&t>t_0=-\frac{1}{\alpha_n}\log
                   {({L_\epsilon}\epsilon)^n}\\
                 \frac{t}{{C_\epsilon}^\frac{1}{n-1}}&t\leq t_0=-\frac{1}{\alpha_n}
                   \log({L_\epsilon}\epsilon)^n
               \end{array}\right.\eqno (3.5)$$
where ${\Lambda_\epsilon},{C_\epsilon},{L_\epsilon}$
will be defined later ( by (3.6), (3.7), (3.10)), which satisfy

i)${L_\epsilon}\rightarrow +\infty$, ${C_\epsilon}\rightarrow +\infty$,
and ${L_\epsilon}\epsilon\rightarrow 0$, as $\epsilon\rightarrow 0$;

ii)${C_\epsilon}-\frac{(n-1)\log(1+c_n{L_\epsilon}^\frac{n}{n-1})+{\Lambda_\epsilon}}{\alpha_n{C_\epsilon}^\frac{1}{n-1}}
=\frac{-\frac{1}{\alpha_n}\log{({L_\epsilon}\epsilon)^n}}{{C_\epsilon}^\frac{1}{n-1}}$;

iii)$\frac{\log{L_\epsilon}}{C_\epsilon^\frac{n}{n-1}}\rightarrow 0$, as
$\epsilon\rightarrow 0$.

Set $u_\epsilon=f_\epsilon({G_p})$, we have $u_\epsilon\in H^{1,n}$ because of
ii).

\vspace{1.5ex}
{\bf Remark 3.1:} When $M=B_1(0)$, the expression of $u_\epsilon$ will become
quite simple:
$$u_\epsilon=\left\{\begin{array}{ll}
                 C_\epsilon-\frac{(n-1)\log(1+c_n|\frac{x}
                    {\epsilon}|^\frac{n}{n-1})+{\Lambda_\epsilon}}{\alpha_n
                   {C_\epsilon}^\frac{1}{n-1}}&|x|\leq L\epsilon\\
                 -\frac{\log{r^n}}{{\alpha_nC_\epsilon}^\frac{1}{n-1}}&
                   |x|> L\epsilon.
               \end{array}\right.$$

We have
$$\begin{array}{lll}
     \int_{A_{t_0}^c}|\nabla u_\epsilon|^ndV_g&=&\int_{A_{t_0}^c}
      \frac{1}{{C_\epsilon}^\frac{n}{n-1}}|\nabla {G_p}|^ndV_g\\[1.5ex]
     &=&\int_{\partial A_{t_0}}\frac{t_0}{{C_\epsilon}^\frac{n}{n-1}}
     |\nabla {G_p}|^{n-2}\frac{\partial {G_p}}
       {\partial n}\\[1.5ex]
     &=&\frac{-\log({L_\epsilon}\epsilon)^n}{\alpha_n{C_\epsilon}^\frac{n}{n-1}}.
  \end{array}$$
and
$$\begin{array}{lll}
     \int_{A_{t_0}}|\nabla u_\epsilon|^ndV_g&=&
       \int_{t_0}^\infty|f'|^n|\nabla {G_p}|^n\int_{\partial A_t}
       \frac{1}{|\nabla {G_p}|}dS_tdt\\[1.5ex]
     &=&\int_{t_0}^\infty|f'|^n\int_{\partial A_t}|\nabla {G_p}|^{n-1}
        dS_tdt\\[1.5ex]
     &=&\int_{t_0}^\infty|f'|^ndt.
  \end{array}$$
Since
$$\begin{array}{llll}
     \int_{t_0}^\infty|f'(t)|^ndt&=&\int_{t_0}^\infty|\frac{c_n\frac{e^{
      \frac{-\alpha_nt}{n-1}}}{\epsilon^\frac{n}{n-1}}}{
      {C_\epsilon}^\frac{1}{n-1}(1+c_n\frac{e^{-\frac{\alpha_nt}{n-1}}}{
      \epsilon^\frac{n}{n-1}})}|^ndt\\[1.5ex]
      &=&\frac{n-1}{\alpha_n{C_\epsilon}^\frac{n}{n-1}}
        \int_0^{c_n{L_\epsilon}^\frac{n}{n-1}}
         \frac{u^{n-1}}{(1+u)^n}du\hspace{4ex}(u=c_n\frac{e^{\frac{-\alpha_nt}{n-1}}}
         {\epsilon^\frac{n}{n-1}})\\[1.5ex]
      &=&\frac{n-1}{\alpha_n{C_\epsilon}^\frac{n}{n-1}}\int_{0}^{c_n{L_\epsilon}^\frac{n}{n-1}}
         \frac{((1+u)-1)^{n-1}}{(1+u)^n}du\\[1.5ex]
      &=&\frac{n-1}{\alpha_n{C_\epsilon}^\frac{n}{n-1}}\sum_{k=0}^{n-2}\frac{C_{n-1}^k(-1)^{
          n-1-k}}{n-k-1}\\[1.5ex]
         &&+\frac{n-1}{\alpha_n{C_\epsilon}^\frac{n}{n-1}}\log(1+c_nL^\frac{n}{n-1})
            +O(\frac{1}{L^\frac{n}{n-1}{C_\epsilon}^\frac{n}{n-1}})\\[1.5ex]
      &=&-\frac{n-1}{\alpha_n{C_\epsilon}^\frac{n}{n-1}}(1+1/2+1/3+\cdots+1/(n-1))\\[1.5ex]
         &&+\frac{n-1}{\alpha_n{C_\epsilon}^\frac{n}{n-1}}\log(1+
         c_n{L_\epsilon}^\frac{n}{n-1})+O(\frac{1}{{L_\epsilon}^\frac{n}{n-1}
         {C_\epsilon}^\frac{n}{n-1}}),
 \end{array}$$
where we use the fact
$$-\sum_{k=0}^{n-2}\frac{C_{n-1}^k(-1)^{
          n-1-k}}{n-k-1}=1+\frac{1}{2}+\cdots+\frac{1}{n-1}.$$
We get
$$\begin{array}{lll}
    \int_M|\nabla u_\epsilon|^ndV_g&=&\frac{1}{\alpha_n{C_\epsilon}^\frac{n}{n-1}
       }\{
       -(n-1)(1+1/2+\cdots+1/(n-1))\\[1.5ex]
       &&+(n-1)\log(1+c_n{L_\epsilon}^\frac{n}{n-1})
       -\log({L_\epsilon}\epsilon)^n+O({L_\epsilon}^{-\frac{n}{n-1}})\}.
   \end{array}$$
Set $\int_M|\nabla u_\epsilon|^ndV_g=1$, we have
$$\begin{array}{lll}
    \alpha_n{C_\epsilon}^\frac{n}{n-1}&=&-(n-1)(1+1/2+\cdots+1/(n-1))\\[1.5ex]
    && +\log\frac{(1+c_n{L_\epsilon}^\frac{n}{n-1})^{n-1}}{{L_\epsilon}^n}-
     \log{\epsilon^n}+O({L_\epsilon}^{-\frac{n}{n-1}})\\[1.5ex]
    &=&-(n-1)(1+1/2+\cdots+1/(n-1))+\\[1.5ex]
    &&\log{\frac{\omega_{n-1}}{n}}-
      \log{\epsilon^n}+O({L_\epsilon}^{-\frac{n}{n-1}}).
  \end{array}\eqno (3.6)$$
By ii), we get
$${\Lambda_\epsilon}=-(n-1)(1+1/2+\cdots+1/(n-1))+O({L_\epsilon}^{-\frac{n}{n-1}}).\eqno (3.7)$$

Next, we compute
$\int_Me^{\alpha_n|u_\epsilon|^\frac{n}{n-1}}dV_g$.

Clearly, $\varphi(t)=|1-t|^\frac{n}{n-1}-\frac{n}{n-1}t
$ is increasing when $0\leq t\leq 1$ and decreasing when
$t\leq 0$, then
$$|1-t|^\frac{n}{n-1}\geq 1-\frac{n}{n-1}t\;\; when\;\; t<1,$$
thus when $t>t_0$, we have
$$\begin{array}{lll}
    \alpha_nf_\epsilon^\frac{n}{n-1}(t)
    &=&\alpha_n{C_\epsilon}^\frac{n}{n-1}|1-\frac{(n-1)\log(1+c_n\frac{e^{-
     \frac{\alpha_n}{n-1}t}}{\epsilon^\frac{n}{n-1}})+{\Lambda_\epsilon}}
     {\alpha_n{C_\epsilon}^\frac{n}{n-1}}|^\frac{n}{n-1}\\[1.5ex]
    &\geq&\alpha_n{C_\epsilon}^\frac{n}{n-1}(1-\frac{n}{n-1}
     \frac{(n-1)\log(1+c_n\frac{e^{-
     \frac{\alpha_n}{n-1}t}}{\epsilon^\frac{n}{n-1}})+{\Lambda_\epsilon}}
     {\alpha_n{C_\epsilon}^\frac{n}{n-1}}).
  \end{array}\eqno(3.8)$$

Applying Lemma 3.1 ,we have
$$\begin{array}{lll}
    \int_{A_{t_0}}e^{\alpha_n|u_\epsilon|^\frac{n}{n-1}}dV_g
     &=&\int_{t_0}^\infty e^{f(t)^\frac{n}{n-1}}\int_{\partial A_t}
     \frac{dS_t}{|\nabla {G_p}|}dt\\[1.5ex]
    &\geq&\int_{t_o}^\infty e^{\alpha_n{C_\epsilon}^\frac{n}{n-1}
     -n\log(1+c_n\frac{e^{\frac{-\alpha_n}{n-1}t}}{
     \epsilon^\frac{n}{n-1}})-\frac{n}{n-1}{\Lambda_\epsilon}
     }\\[1.5ex]
     &&\times e^{-\alpha_nt+\alpha_nS_p}\omega_{n-1}^\frac{n}{n-1}
     (1+O(e^{-\frac{2\alpha_n}{n}t}))\\[1.5ex]
    &=&(1+O(e^{-\frac{2\alpha_n}{n}t_0}))e^{\alpha_n{C_\epsilon}^\frac{n}{n-1}+\alpha_nS_p
     -\frac{n}{n-1}{\Lambda_\epsilon}}\\[1.5ex]
     &&\times\int_{t_0}^\infty\frac{e^{-\alpha_nt}}{(1+c_n\frac{e^{\frac{
     -\alpha_n}{n-1}t}}{
     \epsilon^\frac{n}{n-1}})^n}\omega_{n-1}^\frac{n}{n-1}dt.
 \end{array}$$
Clearly,
$$\begin{array}{llll}
    \int_{t_0}^\infty\frac{e^{-\alpha_nt}}{(1+c_n\frac{e^\frac{
     -\alpha_nt}{n-1}}{\epsilon^\frac{n}{n-1}})^n}\omega_{n-1}^\frac{
     n}{n-1}dt&=&(n-1)\epsilon^n\int_{0}^{c_n{L_\epsilon}^\frac{n}{n-1}}
     \frac{u^{n-2}}{(1+u)^ndu}&(u=c_n\frac{e^{-\frac{\alpha_n}{n-1}t}}{
     \epsilon^\frac{n}{n-1}})\\[1.5ex]
    &=&(n-1)\epsilon^n\int_{0}^{c_n{L_\epsilon}^\frac{n}{n-1}}
     \frac{((u+1)-1)^{n-2}}{(1+u)^n}du\\[1.5ex]
    &=&\epsilon^n(1+O({L_\epsilon}^{-\frac{n}{n-1}})).
  \end{array}$$
Here, we used the fact
$$\sum_{k=0}^m\frac{(-1)^{m-k}}{m-k+1}C_m^k=\frac{1}{m+1}.$$

Then, applying (3.6) and (3.7), we have
$$\begin{array}{lll}
   \int_{A_{t_0}}e^{\alpha_n|u_\epsilon|^\frac{
   n}{n-1}}dV_g&\geq&\frac{\omega_{n-1}}{n}e^{
    \alpha_nS_p+1+1/2+\cdots+1/(n-1)}\\[1.5ex]
    &&\times (1+O(({L_\epsilon}\epsilon)^2))(1+O({L_\epsilon}^{-\frac{n}{n-1}}))\\[1.5ex]
    &=&\frac{\omega_{n-1}}{n}e^{\alpha_nS_p+1+1/2+\cdots+1/(n-1)}+
    \\[1.5ex]
    &&+O(({L_\epsilon}\epsilon)^2
   )+O({L_\epsilon}^{-\frac{n}{n-1}}).
  \end{array}$$
From iii), it follows that $\frac{C_\epsilon}{2}<u_\epsilon|_{A_{t_0}}<
2C_\epsilon$. Hence, we have
$$\int_{A_{t_0}}g_m(u_\epsilon)=O(C_\epsilon^\frac{mn}{n-1}L^2\epsilon^2).$$
Moreover, we have
$$\int_{A_{t_0}^c}(e^{\alpha_n|u_\epsilon|^\frac{
n}{n-1}}-g_m(u_\epsilon))dV_g\geq \mu(M)-\mu(A_{t_0})+\int_{A_{t_0}^c}
\frac{\alpha_n^{m+1}}{(m+1)!}\left|\frac{{G_p}}{{C_\epsilon}^\frac{1}{n-1}}\right|^\frac{n(m+1)}
{n-1}dV_g,$$
thus, we get
$$\begin{array}{l}
    \int_M(e^{\alpha_n|u_\epsilon|^\frac{
        n}{n-1}}-g_m(u_\epsilon))dV_g\\[1.5ex]
        \;\;\;\;\geq\mu(M)+\frac{\omega_{n-1}}{n}e^{
  \alpha_nS_p+1+1/2+\cdots+1/(n-1)}+\\[1.5ex]
  \;\;\;\;\;\;\;\left(\int_{A_{t_0}^c}\frac{\alpha_n^{m+1}}{(m+1)!}|\frac{{G_p}}{{C_\epsilon}^\frac{1}{n-1}}
     |^\frac{n(m+1)}
    {n-1}dV_g+O(C_\epsilon^\frac{mn}{n-1}L^2\epsilon^2)+
  O(({L_\epsilon}\epsilon)^2)+O({L_\epsilon}^{-\frac{n}{n-1}})
  \right)\\[1.5ex]
  \;\;\;\;=\mu(M)+\frac{\omega_{n-1}}{n}e^{
  \alpha_nS_p+1+1/2+\cdots+1/(n-1)}+\\[1.5ex]
  \;\;\;\;\;\;\;{C_\epsilon}^{-\frac{n(m+1)}{(n-1)^2}}
  \left(\int_{A_{t_0}^c}\frac{|\alpha_nG_p^\frac{n}{n-1}|^{m+1}}{(m+1)!}dV_g
   +O({C_\epsilon}^\frac{n(m+1)}{(n-1)^2}({L_\epsilon}\epsilon)^2)
  +O({C_\epsilon}^\frac{n(m+1)}{(n-1)^2}{L_\epsilon}^{-
  \frac{n}{n-1}})\right).
 \end{array}\eqno(3.9)$$
Let
$${L_\epsilon}=(-\log{\epsilon})^{m+1},\eqno (3.10)$$
then as $\epsilon\rightarrow 0$, we have
$${C_\epsilon}^\frac{n(m+1)}{(n-1)^2}{L_\epsilon}^{-\frac{n}{n-1}}\rightarrow 0,\;
  {C_\epsilon}^\frac{n(m+1)}{(n-1)^2}({L_\epsilon}\epsilon)^2\rightarrow 0
  \; and\; \frac{\log{L}}{C_\epsilon^\frac{n}{n-1}}\rightarrow 0.$$
Therefore, i),ii),iii) holds and we can conclude from (3.9) that
$$\int_MF_{1,m}(u_\epsilon)dV_g>\mu(M)+\frac{\omega_{n-1}}{n}e^{
  \alpha_nS_p+1+1/2+\cdots+1/(n-1)}.$$

\end{document}